\documentclass[12pt]{article}
\usepackage{amsfonts,amsmath,amssymb,bbm}
\usepackage{setspace} 
\def\proof{\noindent{\bf Proof:}\hskip10pt}        
\def\QED{\hfill $\Box$}

\font\tenmath=msbm10 scaled 1200
\font\sevenmath=msbm7 scaled 1200
\font\fivemath=msbm5 scaled 1200
\newfam\mathfam \textfont\mathfam=\tenmath
\scriptfont\mathfam=\sevenmath \scriptscriptfont\mathfam=\fivemath

\begin{document}
\def \\ { \cr }
\def\R{\mathbb{R}}
\def \1{1 \mkern -6mu 1} 
\def\N{\mathbb{N}}
\def\E{\mathbb{E}}
\def\P{\mathbb{P}}
\def\Z{\mathbb{Z}}
\def\Q{\mathbb{Q}}
\def\C{\mathbb{C}}
\def\D{\mathbb{D}}
\def\T{\mathbb{T}}
\def \m{^{\rm (m)}}
\def \c{^{\rm (c)}}
\def \+{^{\rm (+)}}
\def \e{{\rm e}}
\def \f{{\mathcal F}}
\def \d{{\tt d}}
\def \r{{\mathcal R}}
\newtheorem{theorem}{Theorem}
\newtheorem{definition}{Definition}
\newtheorem{proposition}{Proposition}
\newtheorem{lemma}{Lemma}
\newtheorem{corollary}{Corollary}
\centerline{\LARGE \bf A limit theorem for trees of alleles}
\vskip 2mm
\centerline{\LARGE \bf  in branching processes}
\vskip 2mm
\centerline{\LARGE \bf  with rare neutral mutations}

\vskip 1cm
\centerline{\Large \bf Jean Bertoin}
\vskip 1cm
\noindent
\centerline{\sl Laboratoire de Probabilit\'es, Universit\'e Pierre et Marie Curie}
\centerline{\sl 175, rue du Chevaleret,   F-75013 Paris, France,}
\centerline{\sl  and DMA, Ecole Normale Sup\'erieure, Paris}
\vskip 15mm

\noindent{\bf Summary. }{\small   We are interested in the genealogical structure of alleles for a Bienaym\'e-Galton-Watson branching process  with neutral mutations (infinite alleles model), in the situation where the initial population is large and the mutation rate small.
We shall establish that for an appropriate regime, the process of the sizes of the allelic sub-families
converges in distribution to a certain continuous state branching process (i.e. a Ji\v{r}ina process) in discrete time. It\^o's excursion theory and the L\'evy-It\^o decomposition of subordinators provide fundamental insights for the results.}
\vskip 3mm
\noindent
 {\bf Key words.}{ \small Weak convergence, branching process, neutral mutations, allelic partition,  L\'evy-It\^o decomposition.} 
 \vskip 5mm
\noindent
{\bf A.M.S. Classification.}  60 J 80, 60 J 05
\vskip 3mm
\noindent{\bf e-mail.} {\tt jean.bertoin@upmc.fr}

\begin{section}{Introduction}
Poisson point processes are a cornerstone of at least two fundamental contributions of Professor Kiyoshi It\^o to Probability Theory,  namely the L\'evy-It\^o decomposition of L\'evy processes (Chapter 1 in \cite{ItoSP}) and It\^o's excursion theory \cite{ItoEx}. The law of rare events, which stresses that 
Poisson variables arise as limiting distributions for the number of successes in a large number of independent trials where each trial has the same small probability of success,
explains their prominent role amongst stochastic processes.
For instance, L\'evy processes can be viewed as weak limits of rescaled random walks and their jumps correspond to rare large steps of the latter. Informally,  the law of rare events thus suggests the Poissonian structure of the jumps of L\'evy processes, and this is indeed the core of the L\'evy-It\^o decomposition. A related but more delicate heuristic also applies to It\^o's description of the excursions of Markov processes, as in discrete time, the succession of the excursions of a Markov chain away from a recurrent point forms an i.i.d. sequence of paths. 

In this paper, we shall  point out the relevance of  this paradigm to study a question motivated by genetics. Recall that Bienaym\'e-Galton-Watson branching processes \cite{AN, Harris, Jagers} model  a population in which at every generation each individual begets according to a fixed offspring distribution and independently of the other individuals, and then dies. Imagine that neutral mutations may happen, so that a child can be either a clone of its parent or a mutant, and the reproduction laws of clones and of mutants are identical. We shall further suppose that each time a mutation occurs, it produces a mutant with a genetic type (allele) which has never been observed before; this setting has been referred to as the infinite alleles model by  Kimura and Crow . 

The allelic partition consists in decomposing the entire population into sub-families of individuals carrying the same allele. 
One important issue in the study of random population models with mutations
(cf. the celebrated sampling formula of Ewens \cite{Ewens} for the Wright-Fisher model) concerns statistics of this allelic partition: what is the probability of observing allelic clusters of certain sizes, how to describe the random genealogical structure connecting these clusters to each other, ...  Our main concern here will be to investigate asymptotics  when the size of the population is large (typically because the number of ancestors is large) and mutations rare. We shall see that, under some mild conditions and for an appropriate regime, a non-degenerate limit exists and is conveniently described in terms of a certain continuous state branching process in discrete time \cite{Jirina}.
It is well-known that continuous state branching processes bear close connexions to certain infinitely divisible distributions; in particular we shall provide a representation of the limiting allelic partitions in terms of Poisson point measures appearing in the L\'evy-It\^o decomposition of the jumps of an underlying L\'evy process.

Let us give a rough idea of the orders of magnitude of the quantities involved. We shall consider a fixed reproduction law with unit mean and finite variance,
and let the Galton-Watson process start from $n$ ancestors having all the same genetic type. It is well-known that if  $n$ generations represent one unit of time and if we rescale the population at each generation by a factor $1/n$, then the rescaled Galton-Watson process
converges in distribution as $n$ tends to infinity to a Feller diffusion.  
We also suppose that neutral mutations affect each child with probability $1/n$. The scaling between population sizes, generations and mutation rates should not come as a surprise since it is precisely the regime of interest for other standard population models, such as the Wright-Fisher model and Kingman coalescent \cite{Kingman}. Recall  that such a critical Galton-Watson process becomes extinct after roughly $n$ generations,  and that the total population is of order $n^2$. So there are only a few mutations at each generation and thus about $n$ different alleles; furthermore the largest allelic sub-families have size of order $n^2$.

Our main result can be described as follows. We use the universal tree ${\mathbb{U}}$, that is the set of finite sequences of integers (including the empty sequence $\varnothing$ that serves as the root of ${\mathbb{U}}$) to record the genealogy of alleles, and define the tree of alleles as a random process ${\mathcal A}$ on ${\mathbb{U}}$, such that the values at vertices are given by the sizes of the corresponding allelic sub-families, with the convention that sizes are ranked in the decreasing order on each sibling. We consider a fixed reproduction law which is critical and has finite variance, and for every integer $n$,  a Galton-Watson process with this reproduction law, started from $n$ ancestors and in which mutations occur at random with rate $1/n$. We write ${\mathcal A}^{(n)}$
for the process on ${\mathbb{U}}$ that describes the corresponding tree of alleles.
Then as $n$ tends to infinity, the rescaled tree of alleles $n^{-2}{\mathcal A}^{(n)}$
converges in the sense of finite dimensional distributions towards a process ${\mathcal A}$ on ${\mathbb{U}}$ with values in $(0,\infty)$. The latter describes the genealogy of a continuous state branching process in discrete time with an inverse Gaussian reproduction law. We stress that its law only depends on the variance of the offspring distribution of the Galton-Watson process, and hence may be viewed as a universal tree of alleles. 
 
 The plan of this paper is as follows. In Section 2, we first present the general setting, stressing the role of the general branching property for the study of Galton-Watson processes with neutral mutations. Then we compute explicitly reproduction laws related to allelic sub-families and point at a connexion with certain downward-skip-free random walks.  Such questions have been addressed from a different point of view in \cite{Be09} to which the present work can be viewed as a complement and a sequel. Section 3 provides some background on continuous state branching processes and convergence of rescaled Galton-Watson processes. The main asymptotic results, namely Proposition \ref{P2} and Theorem \ref{T1} are stated and then proved in Section 4.
\end{section}

 \begin{section}{Galton-Watson processes with neutral mutations}
 \subsection{Basic definitions and branching properties}

 In a Galton-Watson process with neutral mutations, every individual reproduces according to the same distribution  and independently of the other individuals, no matter whether it is a mutant or a clone.  Of course, a clone child of a mutant bears the same allele as its parent. Recall also that we are working in the infinite alleles setting, i.e. the same genetic type cannot be recovered from a cycle of mutations.
 Our basic data are hence provided by a pair of non-negative integer-valued random variables
 $$\xi=(\xi\c,\xi\m)$$
 which describes the number of clone-children and the number of mutant-children of a typical individual. In this paper, we shall mainly be interested in a special situation which appears commonly as a model in population genetics,  namely where mutations affect each child according to a fixed probability and independently of the other children (in other words, the conditional distribution of $\xi\m$ given $\xi\c+\xi\m=\ell$
is binomial with parameter $(\ell,p)$). However the first steps of the analysis can be carried on without difficulties using the general framework.
 We assume throughout this work that
 $$\E(\xi\c)\leq 1\,,$$ i.e. the process of clones is critical or sub-critical
 \footnote{
 Note that this is weaker than assuming that  $\E(\xi\c+\xi\m)\leq 1$ which was required in \cite{Be09}.};
 and we further implicitly exclude the degenerate cases when 
 $\xi\c\equiv 0$, or $\xi\m\equiv 0$.
For every integer $a\geq 1$, we denote by $\P_a$ the law of a Galton-Watson process with neutral mutations, started from $a$ ancestors having the same genetic type and with reproduction law given by that of $\xi=(\xi\c,\xi\m)$.

  The basic {\it branching property} states that for 
every fixed generation, conditionally on the number of individuals at that generation, 
the descents of those individuals are given by independent copies of the initial process,
independently of the preceding generations. It is natural to expect that this
branching property should hold more generally for certain stopping rules, and that this is indeed the case will play an important role in our analysis.
For the sake of simplicity, we shall now present such an extension in a rather informal way, referring to Chauvin \cite{Chauvin} for technical details.

The genealogy of each ancestor is conveniently described by a planar rooted tree, with edges connecting parents to children. More precisely, this requires an additional ordering of the children of each individual, and in this direction we may decide to rank siblings uniformly at random. A line is defined as a family of edges such that every branch from the root (i.e. the ancestor) contains at most one edge in that family. For instance, the edges between parents at generation $k\in\Z_+$ and children at generation $k+1$ form a line. A stopping line should be thought of as a random line such that for every edge in the tree, the event that this edge is part of the line only depends on the marks found on the path from the root to that edge.  
Recall that every edge of the tree corresponds to a pair of individuals  (parent,child), and denote by $C_{\tau}$ the subset of children in the family of edges of some stopping line $\tau$. By removing the edges of $\tau$, we disconnect the genealogical tree into sub-trees whose roots are formed on  the one hand  by the ancestor, and on the other hand by the individuals in $C_{\tau}$. 
The {\it general branching property} then states that conditionally on $C_{\tau}$, the sub-trees rooted at the individuals in $C_{\tau}$ are independent copies of the initial genealogical tree,
and also independent of the initial tree pruned along $\tau$.

We now take into account mutations by assigning marks to the edges between parents and their mutant children. 
Since we are interested by the genealogy of alleles (or equivalently, of mutants), it is convenient  to say that an individual has the $k$-th type if its genotype has been affected by $k$ mutations, that is if its ancestral line comprises exactly $k$ marks. 
Plainly, the family $\tau(k)$ of  edges  connecting a parent of  the $(k-1)$-th type to a mutant child is a stopping line, and the set $C_{\tau(k)}$ coincides with that of
 the mutants of the $k$-th type.  We denote by $T_k$ the total population of individuals of the $k$-th type and by $M_{k}$ the total number of mutants of $k$-th type, agreeing that mutants of the 
$0$-th type are the ancestors (so $M_0=a$, $\P_a$-a.s.). 
The general branching property should make the following statement obvious; we refer the reader to e.g. Chapter Ten in Taib \cite{Taib} for a rigorous argument.

\begin{lemma}\label{L1}
 Under $\P_a$, 
 $$(M_k, k\in\Z_+)$$ 
 is a standard Galton-Watson process 
with reproduction law $\P_1(M_1\in\cdot)$. More generally,
$$((T_k,M_{k+1}), k\in\Z_+)$$
 is a Markov chain 
with transition probabilities 
$$\P_a(T_k=n',M_{k+1}=m'\mid T_{k-1}=n,M_k=m)\,=\,\P_m(T_0=n', M_1=m')\,.$$
\end{lemma}

\noindent{\bf Remark.} We stress the fact that the chain $(T_k, k\in\Z_+)$ of the sizes of sub-populations with given types is not Markov; nonetheless it can be viewed as a {\it hidden} Markov chain. 
In this direction, we also point out that $((T_k,M_k), k\in\Z_+)$ is Markovian,
since the transition probabilities of the chain $((T_{k},M_{k+1}), k\in\Z_+)$ only depend on the second coordinate. Indeed, by a straightforward application of the general branching property, one gets (assuming implicitly that the events on which we condition have positive probability)
\begin{eqnarray*}
& &\P_a(T_{k+1}=n',M_{k+1}=m' \mid M_k=m, T_k=n)\\
&=&\, \P_{m'}(T_0=n') \frac{\P_a(T_{k}=n,M_{k+1}=m' \mid M_k=m)}{ 
\P_a(T_{k}=n \mid M_k=m)}\\
&=&\, \frac{ \P_{m'}(T_0=n')}{ 
 \P_{m}(T_0=n)} \P_m(T_0=n, M_1=m')\,.
\end{eqnarray*}

 Next, observe that  in an infinite alleles model,
the genealogy of individuals naturally induces a genealogy for the alleles in that population. Indeed,
we may identify alleles and mutants, which enables us to use the set of new 
mutants plus a root corresponding to the ancestors of the population (recall that we assume that all the ancestors have the same genetic type) as the set of vertices.
We draw an edge between the root and mutants of the $1$st type, and for every $k\geq 1$ we also draw an edge between 
a mutant of the $k$-th type and a mutant of the $(k+1)$-type if and only if the path
connecting these individuals in the genealogical tree does not contains other mutants.
Hence the set of alleles has a natural structure of  rooted tree. Note that for $k\geq 1$,  $M_k$ corresponds to the number of vertices at 
the $k$-th level
\footnote{For the sake of clarity we shall keep the name {\it generation}
for the distance to the root of individuals in the genealogical tree, and use the name {\it level} when dealing with the structure of alleles.}
 in the tree of alleles. 

Our main goal in this paper is to establish asymptotic features on the genealogy of allelic sub-families,
and in this direction, it will be convenient to view the latter as random processes indexed by the universal tree. More precisely, introduce the set of  finite sequences of positive integers
$${\mathbb{U}}:=\bigcup_{k\in\Z_+}\N^k\,,$$
where $\N=\{1,2,\ldots\}$ and  $\N^0=\{\varnothing\}$. 
Let us briefly recall some standard notation in this setting: 
if  $u=(u_1,\ldots, u_k)$ is vertex  at level  $k\geq 0$ in ${\mathbb{U}}$, then the children of $u$ are $uj:=(u_1,\ldots,u_k,j)$ for $j\in\N$. We also denote by $|u|$ the level of the vertex $u$, with the convention that
the root has level $0$, i.e. $|\varnothing| =0$.
We now take advantage of the natural tree structure of ${\mathbb{U}}$ to record the genealogy of  allelic sub-families together with their sizes.

Given a Galton-Watson process with neutral mutations, we construct  recursively a process ${\mathcal A}=
\left({\mathcal A}_u: u\in {\mathbb{U}}\right)$  as follows.
First, 
${\mathcal A}_{\varnothing}=T_0$
is the size of the sub-population without mutation.
Next, recall that $M_1$ denotes the number of mutants of the first type. 
We enumerate the $M_1$ allelic sub-populations of the first type in the decreasing order of their sizes, with the convention that in the case of ties, sub-populations of the same size are ranked uniformly at random. We denote by
${\mathcal A}_{j}$ the size of the $j$-th allelic sub-populations of the first type, agreeing  that ${\mathcal A}_{j}=0$ if $ j>M_1$. We then complete the construction at all levels by iteration in an obvious way. Specifically, if $ {\mathcal A}_u=0$ for some $u\in{\mathbb{U}}$, then
 $ {\mathcal A}_{uj}=0$ for all $j\in\N$. Otherwise, we enumerate in the decreasing order of their sizes the allelic sub-populations of type $|u|+1$ which descend from the allelic sub-family indexed by the vertex $u$, and then ${\mathcal A}_{uj}$ is the size of this $j$-th sub-family
(as before, in the case of ties, sub-families are ordered uniformly at random, and empty sub-families have size $0$). See Figure 1 for an example. We call the process ${\mathcal A}=
\left({\mathcal A}_u: u\in {\mathbb{U}}\right)$ the {\it tree of alleles}.

\begin{picture}(400,280)(-10,20)

\put (145,40){\makebox(0,0){$\bullet\ $}}

\put (63,120){\makebox(0,0){$\bullet$}}
\put (143,120){\makebox(0,0){$\bullet $}}
\put (223,120){\makebox(0,0){$\clubsuit$}}

\put (143,40){\line(0,1){80}}
\put (143,40){\line(-1,1){80}}
\put (143,40){\line(1,1){80}}

\put (63,121){\line(-1,1){57}}
\put (143,120){\line(-1,1){60}}
\put (143,120){\line(0,1){60}}
\put (223,120){\line(0,1){55}}

\put (3,180){\makebox(0,0){$\heartsuit$}}
\put (83,180){\makebox(0,0){$\bullet $}}
\put (143,180){\makebox(0,0){$\spadesuit$}}
\put (223,180){\makebox(0,0){$\bigcirc$}}

\put (143,180){\line(-1,1){60}}
\put (143,180){\line(0,1){57}}
\put (143,180){\line(1,1){58}}
\put (3,183){\line(0,1){56}}

\put (3,240){\makebox(0,0){$\heartsuit$}}
\put (83,240){\makebox(0,0){$\spadesuit $}}
\put (143,240){\makebox(0,0){$\spadesuit$}}
\put (203,240){\makebox(0,0){$\Diamond$}}

\put (345,40){\circle{20}}
\put (345,120){\circle{20}}
\put (385,120){\circle{20}}
\put (425,120){\circle{20}}
\put (345,200){\circle{20}}
\put (425,200){\circle{20}}

\put (345,40){\makebox(0,0){4}}
\put (425,120){\makebox(0,0){1}}
\put (345,120){\makebox(0,0){3}}
\put (385,120){\makebox(0,0){2}}
\put (345,200){\makebox(0,0){1}}
\put (425,200){\makebox(0,0){1}}

\put (345,50){\line(0,1){60}}
\put (348,48){\line(1,2){32}}
\put (352,46){\line(1,1){67}}
\put (345,130){\line(0,1){60}}
\put (425,130){\line(0,1){60}}

\end{picture}

\vskip 6mm

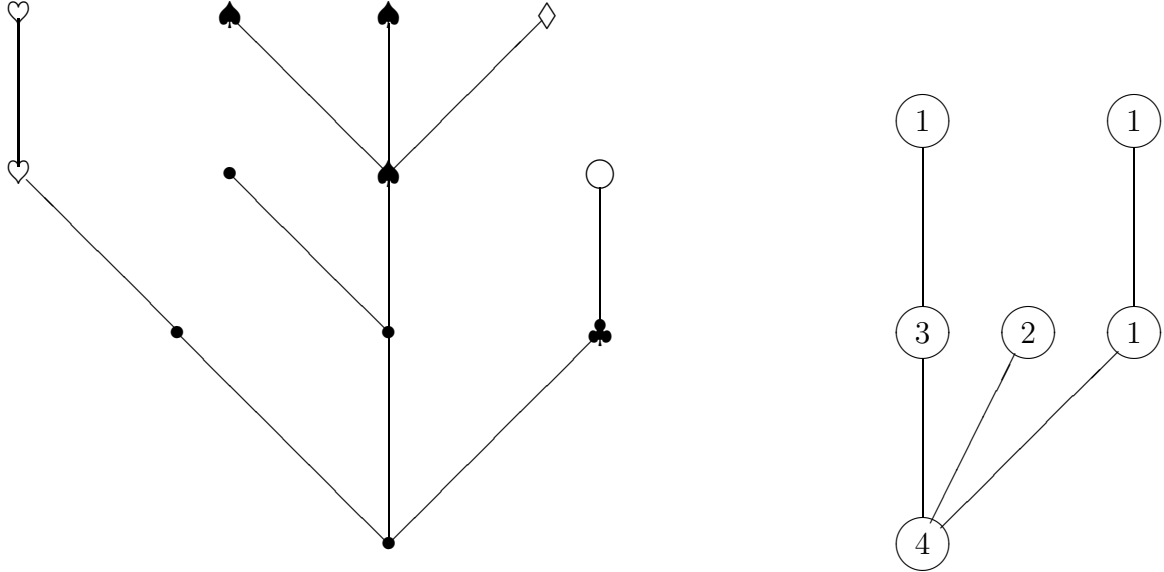
\noindent {\bf Figure 1:  \sl Genealogical tree with mutations (left) and tree of alleles (right).  
The symbols $\bullet, \spadesuit, \heartsuit, \Diamond, \clubsuit, \bigcirc$ represent the different alleles. The labels on the tree of alleles are the sizes of the corresponding allelic sub-families; sub-families with zero size (i.e. which are empty) are omitted.}
\vskip 4mm

It is important to observe that the transition probabilities of the chain 
$((T_k,M_{k+1}), k\in\Z_+)$ in Lemma \ref{L1} depend only on the second coordinate,
and that the latter alone is a Galton-Watson process.  This suggests that the tree of alleles should enjoy some kind of branching property. In order to give a formal statement, is convenient to define first 
the (outer) degree of the tree of alleles ${\mathcal A}$ at some vertex $u\in{\mathbb{U}}$
as
$$ d_u:=\max\{j\geq 1: {\mathcal A}_{uj}>0\}\,,$$
where we agree that $\max \varnothing =0$.
In words,  $d_u$ is the number of allelic sub-populations of type $|u|+1$ which descend from the allelic sub-family indexed by the vertex $u$; in particular $d_{\varnothing}=M_1$.
We shall also need the following notation. Let $\gamma$ be a random variable in
$\N^2$, $d\geq 1$ an integer, and  $\gamma^{(d)}=(\gamma_1, \ldots, \gamma_d)$ where the $\gamma_i$ are  independent copies of $\gamma$. We then denote by $\gamma^{(d\downarrow)}$  the rearrangement  of $\gamma^{(d)}$ in the decreasing order of the first coordinate, with the convention that in the case of ties, the variables $\gamma_i$ with the same first coordinate are ranked uniformly at random.

 The characterization of the probabilistic structure of the tree of alleles that we are now ready to present stems again easily from the general branching property by iteration.  

\begin{lemma}\label{L2}  For every integers $a\geq 1$ and $k\geq 0$,  the tree of alleles fulfills the following properties under $\P_a$ conditionally on $\left(({\mathcal A}_{v}, d_v): |v|\leq k \right)$:

\noindent {\rm (i)} the families of variables
$$\left(({\mathcal A}_{uj},d_{uj}): 1\leq j \leq d_u \right)\,,\qquad\hbox{ $u$ vertex at level $k$ such that }{\mathcal A}_{u}>0\,,$$
are independent, 

\noindent {\rm (ii)}  for each vertex $u$ at level $k$ with ${\mathcal A}_{u}>0$, the $d_u$-tuple 
$\left(({\mathcal A}_{uj},d_{uj}): 1\leq j \leq d_u \right)$
is distributed as $(T_0,M_1)^{(d_u\downarrow)}$ under $\P_1$.
\end{lemma} 

Of course Lemma \ref{L2} of much more informative than the sole Markovian description of the chain $((T_k,M_{k+1}), k\in\Z_+)$ in Lemma \ref{L1} as it retains the information about the genealogy of the allelic sub-families and not merely the sizes of populations of a given type. In this direction,  observe that
$$ T_k=\sum_{|u|=k}  {\mathcal A}_{u}\quad \hbox{and}\quad M_{k+1}=\sum_{|u|=k}d_u\,.$$

  \subsection{Calculation of reproduction laws}
  We shall now determine the transition probabilities that appear in Lemma \ref{L1}.  Essentially, this has been achieved recently in \cite{Be09} using an approach that largely relies on Harris connexion between downward-skip-free random walks and standard Galton-Watson processes, extended to encompass the situation where neutral mutations occur. Here, we shall use a different route, developing calculations that involve generating functions in the case when mutants are supposed to be sterile.

We denote the law of $\xi=(\xi\c,\xi\m)$ by $\pi=(\pi_{k,\ell}: k,\ell\in\Z_+)$, that is
 $$\pi_{k,\ell}:=\P(\xi\c=k, \xi\m=\ell)\,.$$
We also introduce the generating function
 $$g(x,y):=\sum_{k,\ell=0}^{\infty}x^k y^\ell\pi_{k,\ell}=
 \E(x^{\xi\c} y^{\xi\m})\,,\qquad x,y\in[0,1]\,.$$
As we are interested in the joint distribution of the total number of individuals of the $0$-th type and the number of mutants of the $1$-st type, we may imagine a two-type branching process such that clones reproduce independently of each other according to the same distribution $\pi$, while mutants are sterile, i.e. have no progeny a.s.  We write $\varphi$ for the generating function of the total population of $0$-th type and the number of mutants when there is a single ancestor, i.e.
 $$\varphi(x,y):=\E_1(x^{T_0}y^{M_1})\,,\qquad x,y\in[0,1]\,,$$
 so that by the branching property, the generating function of $(T_0,M_1)$ under $\P_a$ is $\varphi^a$.
  The following result is a slight extension of Theorem 1(ii) of \cite{Be09}
 (recall that here we only assume that $\E(\xi\c)\leq 1$ and have an arbitrary number of ancestors, while in \cite{Be09} we worked with a single ancestor and assumed that
 $\E(\xi\c+\xi\m)\leq 1$). It can be viewed as a generalization of the well-known  Dwass formula \cite{Dwass} for the distribution of the total population in standard Galton-Watson processes.

\begin{proposition} \label{P1}
\noindent {\rm (i)} The generating function $\varphi$ is determined by the equation
$$\varphi(x,y)=xg(\varphi(x,y),y)\,,\qquad x,y\in[0,1].$$

\noindent {\rm (ii)} The distribution of $(T_0,M_1)$ is given by
$$\P_a(T_0=n, M_1=\ell)=\frac{a}{n}\pi^{*n}_{n-a,\ell}\,, \qquad n\geq a \geq 1\hbox{ and } \ell\geq 0\,,$$
where $\pi^{*n}$ denotes the $n$-th convolution power of $\pi$ (i.e. $\pi^{*n}$ is the distribution of the sum of $n$ i.i.d. copies of $\xi$).

\end{proposition}

\proof (i) A standard application of the branching property at the first generation gives
\begin{eqnarray*}
\varphi(x,y)&=&\E_1(x^{T_0}y^{M_1})\\
&=&x\sum_{i,j=0}^{\infty}\P(\xi\c=i,\xi\m=j)\E_1(x^{T_0}y^{M_1})^i y^j\\
&=&xg(\varphi(x,y),y)\,.
\end{eqnarray*}
This invites us to consider the equation in the variable $z\in[0,1]$
\begin{equation}\label{E1}
\frac{g(z,y)}{z}=\frac{1}{x}\,,
\end{equation}
where $x,y\in(0,1]$ are fixed. Our assumptions $\E(\xi\c)\leq 1$ and $\xi\c\not\equiv 1$ 
imply that $g(0,y)>0$, and hence $\lim_{z\to 0+}g(z,y)/z=\infty$. On the other hand,
the derivative of $z\to z^{-1}g(z,y)$
is $z\to z^{-2}(z\partial_zg(z,y)-g(z,y))$, and this derivative is strictly negative when
$z>0$ is sufficiently small. This ensures that for each fixed $y\in[0,1]$ and $x>0$  small enough, the equation \eqref{E1} has a unique solution $z=\varphi(x,y)$, and this suffices to determine the law of $(T_0,M_1)$.

 (ii) We shall now derive explicitly the law of $(T_0,M_1)$ under $\P_a$ from its generating function $\varphi^a$ using the classical Lagrange inversion formula. For each fixed $y\in[0,1]$, the function
 $x\to g(x,y)$ is analytic with $g(0,y)\neq 0$. More precisely, we have
  $$g(x,y)=\sum_{k=0}^{\infty}a_k(y) x^k \quad \hbox{with}\quad 
 a_k(y):=\sum_{\ell=0}^\infty  y^\ell\pi_{k,\ell}\,.$$
 According to Lagrange inversion formula (see for instance Section 5.1 in \cite{Wilf}), 
 the $a$-th power of the solution to the equation \eqref{E1} with $y\in[0,1]$ fixed and  $x>0$ sufficiently small, can be expressed in the form
$$\varphi^a(x,y)=\sum_{n=1}^{\infty}\frac{a}{n}\alpha_{n-a}^{*n}\, x^n\,,$$
where $\alpha^{*n}$ stands for the $n$-th convolution power of the finite measure $\alpha=(a_k(y): k\in\Z_+)$.
Observe that the generating function of $\alpha$ is $x\to g(x,y)$, so that of $\alpha^{*n}$
is 
$$x\to g(x,y)^n=\sum_{k=0}^{\infty}x^k\left(\sum_{\ell=0}^\infty  y^\ell\pi^{*n}_{k,\ell}\right)\,.$$ 
Hence we have
$$\alpha_k^{*n}=\sum_{\ell=0}^\infty  y^\ell\pi^{*n}_{k,\ell}\,,$$
and we conclude that
 $$\varphi^a(x,y)=\sum_{n=1}^{\infty}\sum_{\ell=0}^\infty \frac{a}{n} \pi^{*n}_{n-a,\ell}\, x^ny^{\ell}\,,
$$
which completes the proof of (ii).  
 \QED
 
 As generating functions easily yield moments of variables, one immediately deduces from Proposition \ref{P1} simple criteria to decide whether the number of mutant children 
 $M$ is critical, sub-critical, or super-critical, or has a finite second moment.
 
\begin{corollary}\label{C1} {\rm (i)} Suppose that the mean number of clone children is sub-critical, i.e. $\E(\xi\c)<1$.
Then 
$$\E_a(M_1)=a\frac{\E(\xi\m)}{1-\E(\xi\c)}=\E(\xi\m)\E_a(T_0)\,,$$
and in particular
$$\E_1(M_1)   \left\{ \begin{matrix}
<1 &\\ =1 &\\ >1
 &\\
\end{matrix}\right.
\quad \Longleftrightarrow \quad
\E(\xi\c+\xi\m) \left\{ \begin{matrix}
<1 &\\ =1 &\\ >1
 &\\
\end{matrix}\right.
$$
Further
$$\E_1(M_1^2)<\infty \quad \Longleftrightarrow \quad   \E((\xi\c+\xi\m)^2)<\infty\,.$$
\noindent{\rm (ii)}  
If  $\E(\xi\c)=1$, then $\E_1(M_1)=\infty$.
\end{corollary}

\proof Recall that the first moment of an integer-valued variable is given by the left-derivative at $1$ of its generating function.
We get from Proposition \ref{P1}(i)
$$\frac{\partial \varphi}{\partial y}(1,y)=\frac{\partial \varphi}{\partial y}(1,y)
\frac{\partial g}{\partial x}(\varphi(1,y),y)+ \frac{\partial g}{\partial y}(1,y)\,.$$
Since $\varphi(1,1)=1$,  
 this identity forces
$$\frac{\partial \varphi}{\partial y}(1,1)=\E_1(M_1)=\infty$$
when  
 $$\frac{\partial g}{\partial x}(1,1)=\E(\xi\c)=1$$
(recall that $\E(\xi\m)>0$ by assumption), whereas it entails
$$\E_1(M_1)=\frac{\E(\xi\m)}{1-\E(\xi\c)}$$
whenever $\E(\xi\c)<1$.

Observe further that the process of the number of clone children is a branching process with offspring distribution given by the law of
$\xi\c$. In particular, in the sub-critical case $\E(\xi\c)<1$, the total population of clones has a finite expectation
given by  $\E_a(T_0)=a/(1-\E(\xi\c))$. The first equivalence in (i) follows readily. Similar calculations involving the second derivative of generating functions yield the second equivalence in (i). \QED

\subsection{Construction from a random walk}
The starting point of this section is the observation that the transition probabilities of
 the Markov chain $((T_k, M_{k+1}): k\in\Z_+)$ have a simple interpretation in terms of random walks. In this direction, let us first introduce some notation.
We consider a sequence $(\xi_n=(\xi\c_n, \xi\m_n): n\in\N)$ of i.i.d. variables with law
$\pi$,  and then the random walk started from $a\geq 1$ and with steps $\xi\c-1$,
$$S_n\c:=a+\xi\c_1+\cdots+\xi\c_n-n\,,\qquad n\in\Z_+\,.$$
It is convenient to use the (slightly abusive) notation $\P_a$ for  the law of $(S\c_n: n\in\Z_+)$. We also define the first hitting times
$$\varsigma(j):=\inf\{n\in\Z_+: S\c_n=-j\}\,,\qquad j\in\Z_+\,,$$
and
$$\Sigma(j):=\sum_{i=1}^{\varsigma(j)}\xi\m_i\,.$$
We stress that our basic assumption $\E(\xi\c)\leq 1$ ensures that the random walk $S\c$ does not drift to $+\infty$, and hence the passage times $\varsigma(j)$ are finite a.s.
The first identity in next lemma can be viewed as a two-dimensional extension of the well-known result of Otter and Dwass (see e.g. Section 6.2 in \cite{PiSF}) which relates the distribution of the total population in a Galton-Watson process to that of the first hitting time of $0$ of a random walk.

\begin{lemma} \label{L3} The pairs of random variables
$$(\varsigma(0),\Sigma(0))\ \hbox{ and }\ (T_0,M_1)$$
have the same distribution under $\P_a$.
Further, the shifted sequence $(\xi_{\varsigma(0)+j}: j\in\N)$ consists of i.i.d. variables with law $\pi$ and  is independent of $(\varsigma(0), \Sigma(0))$.
\end{lemma}

\proof Introduce for $a=1$ the generating function
$$\tilde \varphi(x,y):=\E_1(x^{\varsigma(0)}y^{\Sigma(0)})\,,\qquad x,y\in[0,1]\,.$$
Because $(S_n: n\in\Z_+)$ is a downwards skip free random walk, an application of the strong Markov property at its first downward passage times shows readily that
for an arbitrary integer $a\geq 1$
$$\E_a(x^{\varsigma(0)}y^{\Sigma(0)})\,=\, \tilde \varphi(x,y)^a\,,\qquad x,y\in[0,1]\,.$$

Now we return to the case $a=1$; 
by conditioning on the first step of the random walk, we get the obvious identity
\begin{eqnarray*}\tilde \varphi(x,y)\,&=&\,\E_1(x^{\varsigma(0)}y^{\Sigma(0)})\\
&=&x\sum_{k,\ell=0}^{\infty}\tilde \varphi(x,y)^k y^{\ell}\pi_{k,\ell}\\
&=&xg(\tilde \varphi(x,y),y)\,,
\end{eqnarray*}
where $g$ denotes the generating function of $\xi=(\xi\c,\xi\m)$. Thus $\tilde \varphi$ solves the equation of Proposition \ref{P1}(i), which establishes our first claim.  As  the hitting time  $\varsigma(0)$ is a stopping time, an application of the strong Markov property then yields the second assertion. \QED

Next, set $\tilde T_0:=\varsigma(0)$, $\tilde M_1:=\Sigma(0)$
 and define  for every $k\in\N$ by an implicit  recurrence
$$
 \tilde T_0+\cdots +\tilde T_k= \varsigma(\tilde M_1+\cdots +\tilde M_k)
 $$
 and 
 $$
\tilde M_1+\cdots+ \tilde M_{k+1}=\Sigma(\tilde T_0+\cdots +\tilde T_k) =\Sigma(\varsigma(\tilde M_1+\cdots +\tilde M_k))\,.$$
Figure 2 below depicts these quantities.

\begin{picture}(300,280)(10,0)

\put (10,160){\line(1,0){460}}
\put (20,5){\line(0,1){250}}

\multiput(49,157)(30,0){14}{$\shortmid$}
\multiput(18,7)(0,30){9}{-}

\put (20,190){\line(1,1){30}}
\put (50,220){\line(1,-1){30}}
\put (80,190){\line(1,-1){30}}
\put (110,159){\line(1,0){30}}
\put (140,160){\line(1,-1){30}}
\put (170,130){\line(1,1){30}}
\put (200,160){\line(1,-1){30}}
\put (230,130){\line(1,-1){30}}
\put (260,100){\line(1,0){30}}
\put (290,100){\line(1,-1){30}}
\put (320,70){\line(1,-1){30}}
\put (350,40){\line(1,1){30}}
\put (380,70){\line(1,-1){30}}
\put (410,40){\line(1,-1){30}}

\multiput(19,97)(10,0){10}{$\cdot$}
\multiput(19,37)(10,0){24}{$\cdot$}
\multiput(19,7)(10,0){34}{$\cdot$}

\put (20,190){\makebox(0,0){$\bullet$}}
\put (50,220){\makebox(0,0){$\bullet$}}
\put (80,190){\makebox(0,0){$\bullet$}}
\put (110,160){\makebox(0,0){$\bullet$}}
\put (140,160){\makebox(0,0){$\bullet$}}
\put (170,130){\makebox(0,0){$\bullet$}}
\put (200,160){\makebox(0,0){$\bullet$}}
\put (230,130){\makebox(0,0){$\bullet$}}
\put (260,100){\makebox(0,0){$\bullet$}}
\put (290,100){\makebox(0,0){$\bullet$}}
\put (320,70){\makebox(0,0){$\bullet$}}
\put (350,40){\makebox(0,0){$\bullet$}}
\put (380,70){\makebox(0,0){$\bullet$}}
\put (410,40){\makebox(0,0){$\bullet$}}
\put (440,10){\makebox(0,0){$\bullet$}}

\put (50,190){\makebox(0,0){$*$}}
\put (110,190){\makebox(0,0){$*$}}
\put (230,220){\makebox(0,0){$*$}}
\put (350,190){\makebox(0,0){$*$}}

\put (110,150){\makebox(0,0){$\varsigma(0)$}}
\put (170,150){\makebox(0,0){$\varsigma(1)$}}
\put (260,150){\makebox(0,0){$\varsigma(2)$}}
\put (320,150){\makebox(0,0){$\varsigma(3)$}}
\put (350,150){\makebox(0,0){$\varsigma(4)$}}
\put (440,150){\makebox(0,0){$\varsigma(5)$}}

\put (5,100){\makebox(0,0){$-2$}}
\put (204,100){\vector(1,0){55}}
\put (168,100){\vector(-1,0){55}}
\put (190,100){\makebox(0,0){$\tilde T_1$}}

\put (5,40){\makebox(0,0){$-4$}}
\put (325,40){\vector(1,0){25}}
\put (280,40){\vector(-1,0){25}}
\put (300,40){\makebox(0,0){$\tilde T_2$}}

\put (5,10){\makebox(0,0){$-5$}}
\put (415,10){\vector(1,0){25}}
\put (373,10){\vector(-1,0){20}}
\put (395,10){\makebox(0,0){$\tilde T_3$}}

\end{picture}

\vskip20pt 
\noindent 
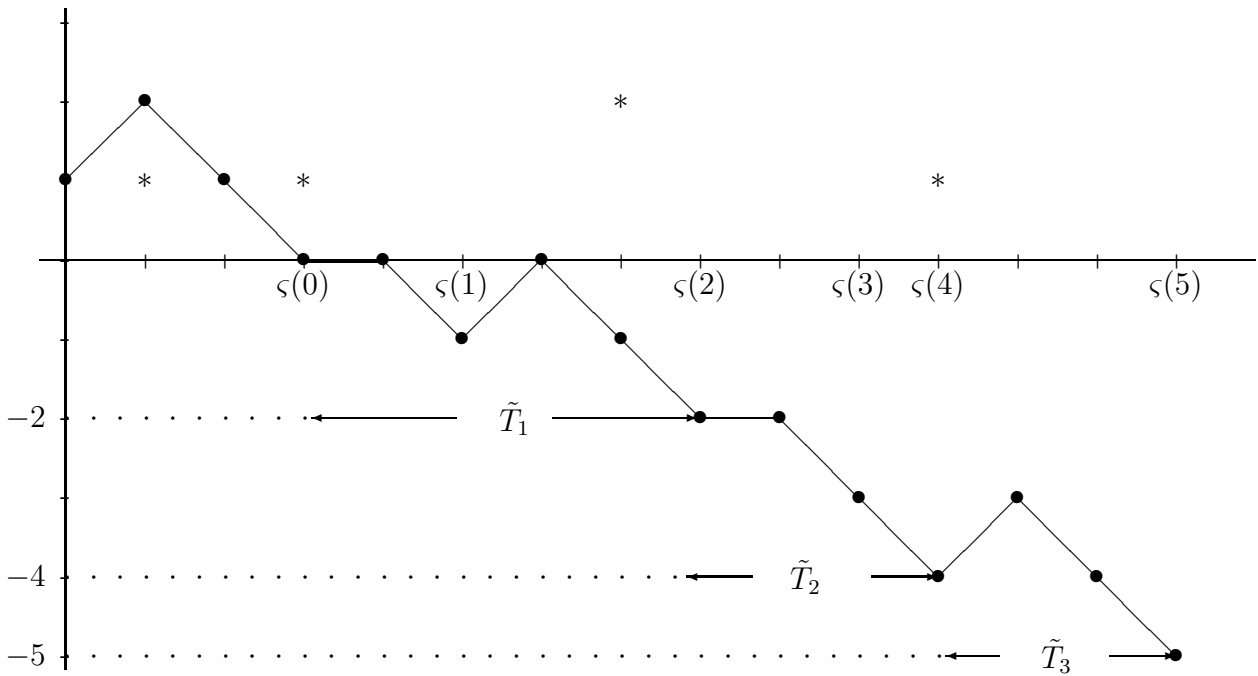
{\bf Figure 2:  \sl The graph of the random walk $S\c$; the $*$ represent the non-zero values of the variables $\xi\m$.  Here $\tilde M_1=2$, $\tilde M_2=2$, $\tilde M_3=1$ and $\tilde M_4=0$.}
\vskip 1cm

\begin{corollary}\label{C4} For every $a\geq 1$, the chains $((T_k,M_{k+1}): k\in\Z_+)$ and $((\tilde T_k,\tilde M_{k+1}): k\in\Z_+)$ have the same distribution under $\P_a$.
\end{corollary}

\proof  It is immediately checked by induction that each $ \tau_k:=\tilde T_0+\cdots +\tilde T_k$ is a stopping time in the natural filtration $({\mathcal G}(n))_{n\in\N}$  generated by the i.i.d. sequence $(\xi_n: n\in\N)$, and that
$M_{k+1}$ is ${\mathcal G}( \tau_k)$-measurable. By an application of the strong Markov property, we get that $((\tilde T_k,\tilde M_{k+1}): k\in\Z_+)$ is a homogeneous Markov chain. 
More precisely, the conditional distribution of $(\tilde T_k,\tilde M_{k+1})$ given 
$\tilde T_{k-1}=t$ and $\tilde M_k=m$ is that of $(\varsigma(0),\Sigma(0))$ under $\P_m$.
Combining these observation with Lemmas \ref{L1} and \ref{L3} completes the proof. 
\QED

More generally, we can apply Lemma \ref{L3} to construct from the i.i.d.  variables
$\xi_n$ a random process ${\mathcal A}'$ indexed by ${\mathbb{U}}$  with the same distribution as  tree of alleles ${\mathcal A}$, by making use of  the characterization of the law of the latter in  Lemma \ref{L2}.
To start with, the process  ${\mathcal A}'$ fulfills the following two requirements.
First, if $ {\mathcal A}'_u=0$ for some $u\in{\mathbb{U}}$, then  $ {\mathcal A}'_{uj}=0$ for all $j\in\N$.
Second,  for every vertex $u\in{\mathbb{U}}$ such that ${\mathcal A}'_u>0$,  the (outer) degree of  $ {\mathcal A}'$ at $u$, 
$$d'_u:=\#\{j\in\N: {\mathcal A}'_{uj}>0\}\,,$$
 is a finite number and ${\mathcal A}'_{uj}>0$ if and only if $j\leq d'_u$.
We set  ${\mathcal A}'_{\varnothing}=\varsigma(0)$ and $d'_{\varnothing}=\Sigma(0)$.
 Next, consider the increments
$$\lambda(j):=\varsigma(j)-\varsigma(j-1) \ \hbox{ and } \ \delta(j)=\Sigma(j)-\Sigma(j-1))\,,\qquad j\geq 1\,,$$
For vertices at the first level, $(({\mathcal A}'_j, d'_j): 1\leq j \leq d'_{\varnothing})$ is given by
the rearrangement of the sequence $((\lambda(j),\delta(j)): 1\leq j \leq d'_{\varnothing})$ in the  decreasing order the first coordinate $\lambda(j)$ (with the usual convention in case of ties). We may then continue with vertices of the next levels by an iteration which should be obvious (but which would also be quite intricate to state explicitly). Figure 3 below may help visualizing the construction.

\begin{picture}(400,330)(-100,0)

\put (145,40){\circle{30}}

\put (63,120){\circle{30}}
\put (223,120){\circle{30}}

\put (131,50){\line(-1,1){58}}
\put (156,51){\line(1,1){56}}
\put (53,131){\line(-1,1){53}}
\put (75,130){\line(1,1){58}}
\put (145,215){\line(0,1){49}}

\put (145,200){\circle{30}}
\put (-5,200){\circle{30}}
\put (145,280){\circle{30}}

\put (145,40){\makebox(0,0){$\varsigma(0)$}}
\put (63,120){\makebox(0,0){$\lambda(2)$}}
\put (223,120){\makebox(0,0){$\lambda(1)$}}
\put (145,200){\makebox(0,0){$\lambda(4)$}}
\put (-5,200){\makebox(0,0){$\lambda(3)$}}
\put (145,280){\makebox(0,0){$\lambda(5)$}}

\end{picture}

\vskip 4mm

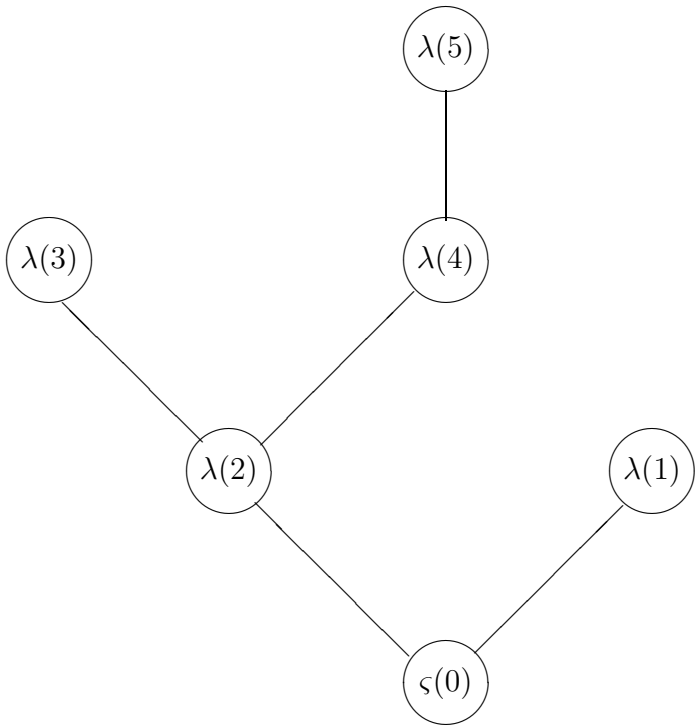
\noindent {\bf Figure 3: \sl Tree of alleles constructed from the random walk $S\c$ and the variables $\xi\m$ of Figure 2. The labels on the vertices are the lengths of the excursions of $S\c$ above its current minimum, they correspond to the sizes of the allelic sub-families (again sub-families with size $0$ are omitted).}
\vskip 4mm

\end{section}

\begin{section}{Background on continuous state branching processes}
Before describing our main limit  results for trees of alleles, we need to develop some basic material about limits of rescaled Galton-Watson processes.
The L\'evy-It\^o decomposition of subordinators plays a crucial role for the representation of the genealogical structure of the continuous state limits.

We start with the classical convergence to Feller diffusions \cite{Feller, Jirina}, i.e. the solutions
$(X(x,t), t\geq 0)$ to stochastic differential equations of the type
\begin{equation}\label{E5}
X(x,t)=x+ \int_0^t \sigma \sqrt{X(x,s)}{\rm d}B_s+b\int_0^t X(x,s){\rm d}s\,,\qquad t\geq 0\,,
\end{equation}
where  $x\geq 0$ is the initial value, $b\in\R$ and $\sigma^2>0$ are parameters, and
 $(B_t: t\geq 0)$ denotes a standard Brownian motion. 
For every $n\in\N$, consider  a Galton-Watson process $(Z^{(n)}_k: k\in\Z_+)$ which starts from 
$Z^{(n)}_0=a(n)$ ancestors and has reproduction law $\rho^{(n)}$, where $\rho^{(n)}$ is some probability measure on $\Z_+$ and $a(n)$ a positive integer.
Write 
$${\rm m}(\rho^{(n)}):=\sum_{i=0}^{\infty}i\rho^{(n)}_i\ \hbox{ and }
{\rm var}(\rho^{(n)}):= \sum_{i=0}^{\infty}(i-{\rm m}(\rho^{(n)}))^2\rho^{(n)}_i$$
for the first moment and the variance of $\rho^{(n)}$.
In the situation where
\begin{equation} \label{E13}
a(n)\sim nx\ ,\ {\rm m}(\rho^{(n)})-1\sim bn^{-1}\ \hbox{ and } 
{\rm var}(\rho^{(n)}) \sim \sigma^2\qquad \hbox{as }n\to\infty
\end{equation}
 for some $x\in(0,\infty)$, $b\in\R$ and $\sigma^2>0$,
it is well known that
\begin{equation}\label{E3}
(n^{-1}Z^{(n)}_{\lfloor nt\rfloor}:t\geq 0) \ \Longrightarrow\ 
(X(x,t): t\geq 0)
\end{equation}
where the notation $ \Rightarrow$ refers to convergence in distribution as $n\to\infty$
and $X(x,t)$ is the Feller diffusion specified by \eqref{E5}.

We next turn our attention to the simpler situation where one only rescales the number of individuals and uses the generations as a discrete time parameter. For the sake of clarity, we shall deal with a framework that is slightly less general than it could be.
We denote the tail distribution of $\rho^{(n)}$ by $\bar\rho^{(n)}(y):=\rho^{(n)}((y,\infty))$ for $y>0$
and now assume  that
\begin{equation}\label{E6}
\lim_{n\to\infty} n^{-1}a(n)=x\ \hbox{ and }\ 
\lim_{n\to\infty} n\bar \rho^{(n)}(ny)= \bar\nu (y)\quad \hbox{in }L_{\rm loc}^1([0,\infty),{\rm d}y)\,,
\end{equation}
where $\bar\nu$ is some locally integrable non-increasing function on $[0,\infty)$ with $\bar \nu(\infty)=0$. We may thus think of $\bar \nu$ as the tail of a Radon measure $\nu$
on $(0,\infty)$ with $\int(1\wedge y)\nu({\rm d}y)<\infty$; $\nu$ will be often referred
 to as a L\'evy measure.
Our assumptions ensure that
\begin{equation}\label{E17}
n^{-1}Z^{(n)}_1\,  \Longrightarrow\  \, Z_1\,,
\end{equation}
where $Z_1$ is a random variable with values in $[0,\infty)$ which is infinitely divisible. 
Indeed, we  have for any $q>0$
\begin{eqnarray*}
\E(\exp(-q n^{-1} Z^{(n)}_1))&=&\left(  1-\int_{[0,\infty)}(1-\e^{-qy/n})\rho_n({\rm d}y)\right) ^{a(n)}\\
&=&\left(  1-\frac{q}{n}\int_{0}^{\infty}\e^{-qy/n}\bar\rho_n(y){\rm d}y\right) ^{a(n)}\\
&=&\left(  1-q\int_{0}^{\infty}\e^{-qy}\bar\rho_n(ny){\rm d}y\right) ^{a(n)}\,,\\
\end{eqnarray*}
and \eqref{E6} ensures that the latter quantity converges as $n\to\infty$ towards the Laplace transform of an infinitely divisible variable 
$$\E(\exp(-qZ_1))=\exp( -x \kappa(q))\,,$$
where the cumulant $\kappa$ is given by the L\'evy-Khintchine formula
\begin{equation}\label{E14}
\kappa(q)=\int_{(0,\infty)}(1-\e^{-qy})\nu({\rm d}y)\,.
\end{equation}
We underline the fact that the drift coefficient is $0$; this will play an important role in the sequel. 
An application of the Markov property now shows that more generally 
\begin{equation}\label{E8}
(n^{-1}Z^{(n)}_k: k\in\Z_+) \ \Longrightarrow\ (Z_k: k\in\Z_+)
\end{equation}
where $(Z_k: k\in\Z_+)$ is a Markov chain with values in $\R_+$, started from 
$Z_0=x$ and whose transition probabilities are characterized as follows: for every 
$k\in\Z_+$ and $q,y\geq 0$,
$$\E(\e^{-qZ_{k+1}}\mid Z_k=y)\,=\, \exp(-y\kappa(q))\,.$$
One refers to the limiting chain $Z$ as a (discrete time) continuous state branching process, in short CSBP, with reproduction measure $\nu$ and started from $x$.

It is interesting to recast the preceding convergence in the framework of the law of rare events. In this direction, recall that 
the L\'evy-It\^o decomposition of the infinitely divisible variable $Z_1$ reads
\begin{equation} \label{E7}
Z_1=\sum_{i=1}^{\infty}{\tt a}_i\,,
\end{equation}
where  ${\tt a}_1\geq {\tt a}_2\geq \ldots$ are the atoms ranked in the decreasing order of a Poisson random measure on $(0,\infty)$ with intensity $x\nu$, with the convention 
that atoms are repeated according to their multiplicity and that
when the Poisson random measure is finite (which occurs if and only if $\nu((0,\infty))<\infty$), then ${\tt a}_i=0$ whenever the index $i$ exceeds the total mass of the Poisson measure. 
Consider for every $n\in\N$ a family $(\xi^{(n)}_i: 1\leq i \leq a(n))$ of i.i.d. variables with law $\rho^{(n)}$; we should think of $\xi^{(n)}_i$ as the number of children of the $i$-th ancestor in the Galton-Watson process $Z^{(n)}$. Denote by 
${\tt a}^{(n)}_1\geq {\tt a}^{(n)}_2\geq \ldots\geq {\tt a}^{(n)}_{a(n)}$ the decreasing reordering of the rescaled variables $(n^{-1}\xi^{(n)}_i: 1\leq i \leq a(n))$.
In the regime \eqref{E6}, the  law of rare events for null arrays (e.g. Theorem 14.18 in \cite{Kallenberg}) ensures that
\begin{equation}\label{E9}
 ({\tt a}^{(n)}_1,{\tt a}^{(n)}_2, \ldots, {\tt a}^{(n)}_{a(n)})  \Longrightarrow\ 
({\tt a}_1, {\tt a}_2, \ldots)\,,
\end{equation}
in the sense of finite dimensional distributions. Note also that \eqref{E17} can be re-written in this setting as 
$$\sum_{i=1}^{a(n)} {\tt a}^{(n)}_i \ \Longrightarrow\, \sum_{i\in\N} {\tt a}_i\,;
$$
however the latter does not follow from \eqref{E9}.

This invites us to describe the convergence of rescaled Galton-Watson processes to (discrete time) CSBP from another point of view
that takes into account the genealogy, and not merely the total sizes of populations at
given generations. In this direction, we use a representation of the latter as random processes indexed by the universal tree ${\mathbb{U}}$. For simplicity, suppose for a while that the L\'evy measure $\nu$ is infinite, so a Poisson random measure with intensity $c\nu$ with $c>0$ has infinitely many atoms a.s.
Recall from the L\'evy-It\^o decomposition \eqref{E7} that almost all the individuals at the first generation in a CSBP descend from only countably many ancestors (we stress again that we are dealing with cumulants $\kappa$ with no drift component), and plainly the same feature holds for the subsequent generations.
Roughly speaking,  vertices $u\in{\mathbb{U}}$ at level $|u|=k\geq 1$ represent the sizes  of the sub-populations at generation $k$ in the CSBP which descent from the same parent at generation $k-1$.
We construct a random process $({\mathcal Z}_u: u\in{\mathbb{U}})$ 
related to the CSBP $Z$, where 
${\mathcal Z}_{uj}$ is the size of the $j$-th largest sub-population at generation $|u|+1$ 
which descents from a parent in the sub-population represented by $u$. We stress the process ${\mathcal Z}$ is by definition non-increasing on each sibling, i.e. the map $j\to {\mathcal Z}_{uj}$ is non-increasing on $\N$ for every $u\in{\mathbb{U}}$.
More precisely, conditionally on ${\mathcal Z}_u=z$, the L\'evy-It\^o decomposition \eqref{E7} suggests that 
 $({\mathcal Z}_{uj}: j\in\N)$ should be given by the sequence of the atoms of a Poisson random measure on $(0,\infty)$ with intensity $z\nu$, where atoms are repeated according to their multiplicity and ranked in the decreasing order. 
 We make the construction formal in the following definition.

\begin{definition}  Fix $x>0$ and $\nu$ a measure on  $(0,\infty)$ with $\int(1\wedge y)\nu({\rm d}y)<\infty$.
A tree-indexed CSBP with reproduction measure $\nu$ and initial population of size $x$
is a process $({\mathcal Z}_{u}: u\in{\mathbb{U}})$ with values in $\R_+$ and indexed by the universal tree,  whose distribution
is characterized by induction on the levels  as follows:

\noindent {\rm (i)} ${\mathcal Z}_{\varnothing}=x$ a.s.;

\noindent {\rm (ii)} for every $k\in\Z_+$, conditionally on $({\mathcal Z}_v: v\in{\mathbb{U}}, |v|\leq k)$, the sequences $({\mathcal Z}_{uj})_{j\in\N}$ for the vertices $u\in{\mathbb{U}}$ at generation $|u|=k$ are independent, and each sequence $({\mathcal Z}_{uj})_{j\in\N}$ is distributed as the family of the atoms of a Poisson random measure on $(0,\infty)$ with intensity ${\mathcal Z}_u \nu$, where atoms are repeated according to their multiplicity,  ranked in the decreasing order, and completed by an infinite sequence of $0$ if the Poisson measure is finite.
\end{definition}

It should be plain that if ${\mathcal Z}$ is a tree-indexed CSBP with reproduction measure $\nu$ and initial population of size $x$, then $\left(
\sum_{|u|=k} {\mathcal Z}_u : k\in\Z_+\right)$ is a CSBP with reproduction measure $\nu$ started from $x$.
We also point out that for every integer $n$, we can represent similarly the genealogy for the Galton-Watson process $Z^{(n)}$ as a process ${\mathcal Z}^{(n)}$ indexed by the universal tree ${\mathbb{U}}$,
and one can check that under the regime \eqref{E6}, the following extension of \eqref{E8} holds:
$$n^{-1}{\mathcal Z}^{(n)}\ \Longrightarrow\ {\mathcal Z}$$
in the sense of finite dimensional distributions. 
This should be viewed as a variation of the law of rare events \eqref{E9}; the easy proof is left to the interested reader.

We now conclude this section by underlying the connexion between discrete time CSBP and subordinators (i.e. L\'evy processes with values in $\R_+$). Consider a subordinator
 $\tau=(\tau_t: t\geq 0)$ with  no drift and L\'evy measure $\nu$. Its cumulant $\kappa$
 is given by the L\'evy-Khintchine formula \eqref{E14} and we have
$$\E(\e^{-q \tau_t})=\exp(-t\kappa(q))\,,\qquad \hbox{for all }q,t\geq 0\,.$$
Fix $x>0$ and define  a sequence $(\zeta_k: k\in\Z_+)$ by implicit iteration as follows:
$$\zeta_0=x\,,\  \zeta_1= \tau_x\,,\  \zeta_1+\zeta_2= \tau_{x+\zeta_1}
\,,\ \ldots\,, \ \zeta_1+\cdots + \zeta_{k+1}=\tau_{x+\zeta_1+\cdots+\zeta_k}\,.$$
Observe by an easy induction that the random times $x+\zeta_1+\cdots+\zeta_k$
are stopping times in the natural filtration of $\tau$, so that the strong Markov property 
can be applied. It is then immediate to check that  $(\zeta_k: k\in\Z_+)$ is a CSBP with reproduction measure $\nu$ and initial population of size $x$. 

More generally, the  tree-indexed CSBP ${\mathcal Z}$ can be constructed from  the subordinator $\tau$ by making full use of the L\'evy-It\^o decomposition. Specifically, 
we know from the latter that the Stieltjes measure $d\tau$ on the random interval 
$$I_k:=(x+\zeta_1+\cdots+\zeta_{k-1}, x+\zeta_1+\cdots+\zeta_k]$$ 
is purely atomic, and conditionally on $|I_k|$, the sequence of the atomic masses has the same distribution as the family of the atoms in a Poisson point measure on $(0,\infty)$ with intensity  $|I_k|\nu$. These atoms should be viewed as the sizes of sub-families at level $k$, so it remains to identify the siblings and rank atoms corresponding to a same sibling in the decreasing order. This is straightforward for the first levels but becomes increasingly intricate for larger levels. 
Specifically, we let ${\mathcal Z}_{\varnothing}=x$ and
declare that $({\mathcal Z}_j: j\in\N)$ is given by the sequence of the jumps of $\tau$ on $(0,x]$ ranked in the decreasing order. Next  $({\mathcal Z}_{1j}: j\in\N)$ corresponds to the ranked sequence of the jumps of $\tau$ on the interval
$(x, x+\tau_{{\mathcal Z}_1}]$, $({\mathcal Z}_{2j}: j\in\N)$ to those on the interval 
$(x+\tau_{{\mathcal Z}_1}, x+\tau_{{\mathcal Z}_1+{\mathcal Z}_2}]$ and so on. The algorithm may be thought of as a variant of the breadth first search in which each sibling is ordered according to the size of its progeny.

\end{section}

\begin{section}{Asymptotic for rare mutations}

This section contains our main results on limits of trees of alleles; we shall first present and discuss the general framework, then state the results, and finally prove the latter.

\subsection{Framework and main results}
We consider a fixed probability measure
$\pi^{(+)}$ on $\Z_+$ which serves as reproduction law for a standard Galton-Watson process denoted by $Z\+$. We assume that $Z^{(+)}$ is critical, i.e.
$$
\sum_{i=0}^{\infty} i\pi^{(+)}_i=1\,,
$$
and has a finite variance
$$\sum_{i=0}^{\infty}(i-1)^2\pi^{(+)}_i=\sigma^2 <\infty\,.$$
Further, we suppose that  mutations affect each child according to a fixed probability $p\in(0,1)$ and independently of the other children. That is to say that the probability measure $\pi$ on $\Z_+\times \Z_+$ which gives the law of the number of clone children and the number of mutant children of a typical individual is given by
$$\pi_{k,\ell}=\pi^{(+)}_{k+\ell}  \left( \begin{matrix}
k+\ell \\ k \\ \end{matrix}\right) (1-p)^k p^{\ell}\,.
$$
We will use the notation $\P^{p}_{a}$ for the probability measure under which the
Galton-Watson process $Z\+$ has $a$ ancestors and the mutation rate is $p$,
and ${\mathcal L}\left (\cdot, \P^{p}_{a}\right)$ will then refer to the distribution of a random variable or a process under $\P^{p}_{a}$

We are interested in the situation where the mutation rate $p=p(n)$ is small and the number of ancestors $a=a(n)$ large when the parameter $n$ goes to infinity. 
Specifically,  we consider the regime
\begin{equation}\label{E11}
a(n)\sim nx \ \hbox{ and }\ p(n)\sim c n^{-1}\,,
\end{equation}
where $c,x$ are some positive constants. Let us start by mentioning some results of convergence 
in distribution for Galton-Watson processes in this setting.

First, we know from \eqref{E3} that the Galton-Watson process $Z\+$ properly rescaled converges to a Feller diffusion  on $\R_+$; specifically
\begin{equation}\label{E10}
{\mathcal L}\left((n^{-1}Z^{(+)}_{\lfloor nt\rfloor}:t\geq 0), \P^{p(n)}_{a(n)}\right) \ \Longrightarrow\ 
(X\+_t: t\geq 0)\,,
\end{equation}
where $(X\+_t: t\geq 0)$ solves the SDE \eqref{E5} for the parameter $b=0$. 
In the same direction, the marginal law of $\xi\c$ under $\P^{p(n)}$ has first moment $1-p(n)$ and variance close to $\sigma^2$ when $n$ is large.
Hence,  if $Z\c$ denotes the Galton-Watson process of clones (i.e. we only consider individuals of the $0$-th type), then 
\begin{equation}\label{E18}
{\mathcal L}\left ((n^{-1}Z\c_{\lfloor nt\rfloor}:t\geq 0), \P^{p(n)}_{a(n)}\right) \ \Longrightarrow\ 
(X\c_t: t\geq 0)\qquad \hbox{as $n\to\infty$},
\end{equation}
  where 
$(X\c_t: t\geq 0)$ is another Feller diffusion
solution to the SDE \eqref{E5}
for the parameter $b=-c$.

On the other hand, recall from Lemma \ref{L1} and  Corollary \ref{C1}(i) that the process of the number of mutants of given types $(M_k: k\in\Z_+)$  is a critical Galton-Watson process with finite variance. In view of the classical limit theorem stated as \eqref{E3} in Section 3, one might suspect that the rescaled process $(n^{-1}M_{\lfloor nt\rfloor}: t\geq 0)$ could converge to some Feller diffusion. However this is not the case; indeed an easy calculation shows that the variance of the reproduction law of $M_{\cdot}$ under $\P^{p(n)}$ is of order $n$, and thus the requirement \eqref{E13} fails. 
Nonetheless one can deduce from a few lines of calculations based on Proposition \ref{P1} that the condition \eqref{E6} is fulfilled by the reproduction law of $M$ under $\P^{(n)}_{a(n)}$, and hence 
$${\mathcal L}\left ((n^{-1}M_{k}): k\in\Z_+), \P^{p(n)}_{a(n)}\right)$$
converges weakly when $n\to\infty$ towards the law of some discrete time CSBP started from $a$. We do not give a formal statement as the forthcoming Proposition \ref{P2} is a stronger result.

The asymptotics \eqref{E10} and \eqref{E18}  point to the fact that in the regime \eqref{E11}, the total size of the population of the Galton-Watson process should be rescaled by a factor $n^{-2}$, and in particular 
the asymptotic behavior of the number $T_0=\sum_{k=0}^{\infty} Z\c_k$ of individuals of $0$-th type is given by
$${\mathcal L}\left (n^{-2}T_0, \P^{p(n)}_{a(n)}\right)   \ \Longrightarrow\ 
\int_0^{\infty}X\c_t {\rm d}t\,.$$
More generally, we have the following joint convergence in distribution for the rescaled process of the sizes of sub-populations and the number of mutants of a given type.

\begin{proposition}\label{P2} In the regime \eqref{E11}, we have
$${\mathcal L}\left (((n^{-2}T_k, n^{-1}M_{k+1}): k\in\Z_+), \P^{p(n)}_{a(n)}\right)   \ \Longrightarrow\ 
\left((Z_{k+1},cZ_{k+1}): k\in\Z_+\right)$$
where $(Z_k: k\in\Z_+)$ is a CSBP with reproduction measure 
$$\nu({\rm d}y) = \frac{c}{\sqrt{2\pi \sigma^2 y^3}}\exp\left(-\frac{c^2 y}{2\sigma^2}\right){\rm d}y\,,
\qquad y>0\,,$$
 and initial population of size $x/c$.
\end{proposition}

The L\'evy-It\^o decomposition now suggests that conditionally on $n^{-2}T_k\sim y$, the sequence of the sizes of the sub-populations carrying a same allele of the $(k+1)$-type and normalized by a factor $n^{-2}$ 
should converge in distribution to the sequence of the atoms of a Poisson random measure  on $\R_+$ with intensity specified in Proposition \ref{P2}. Recall also that $d_u$ denotes the outer degree at the vertex $u\in{\mathbb{U}}$ in the tree of alleles, and observe from Lemma \ref{L2} that for a Galton-Watson process with neutral mutations, the process $(({\mathcal A}_u,d_u): u\in{\mathbb{U}})$  has a simpler Markovian structure than $({\mathcal A}_u: u\in{\mathbb{U}})$ alone. This leads us to our main limit theorem for the tree of alleles. 

\begin{theorem}\label{T1} In the regime \eqref{E11}, the rescaled tree of alleles $n^{-2}{\mathcal A}$ under $\P_{a(n)}^{p(n)}$ converges in the sense of finite dimensional distributions to the tree indexed CSBP $({\mathcal Z}_u : u\in{\mathbb{U}})$ with reproduction measure $\nu$ given in Proposition \ref{P2} and  random initial population with inverse Gaussian distribution:
$$\frac{\P({\mathcal Z}_{\varnothing}\in {\rm d}y)}{{\rm d}y}= 
\frac{x}{\sqrt{2\pi \sigma^2 y^3}}\exp\left(-\frac{(cy-x)^2}{2\sigma^2y}\right)\,,\qquad y>0\,.$$

More precisely, if we also take into account the outer degrees, then we have the joint convergence in the sense of finite dimensional distributions:
$${\mathcal L}\left (((n^{-2}{\mathcal A}_u, n^{-1}d_u): u\in{\mathbb{U}}), \P^{p(n)}_{a(n)}\right)   \ \Longrightarrow\ 
\left({\mathcal Z}_u , c{\mathcal Z}_u): u\in{\mathbb{U}}\right)\,.$$

\end{theorem}

\subsection{Proofs}
Let us first present informally some intuitions for the proofs, which rely on the connexion with random walks in Section 2.3. Roughly speaking,  we shall observe that in the regime \eqref{E11}, the random walk $S\c$ suitably rescaled converges to a Brownian motion with negative drift. As the lengths of the excursions of $S\c$ above its current minimum correspond to the sizes of sub-populations with the same allele, this suggests that in the limit, the lengths of the excursions of a Brownian motion with drift above its current minimum should describe the limit of rescaled sub-populations. According to It\^o's excursion theory, these lengths  can be described in terms of a Poisson point process. The comparison with the construction of the tree indexed CSBP presented in Section 3.2 should then make  Theorem \ref{T1} more intuitive.

The proofs of Proposition \ref{P2}  and Theorem \ref{T1} both  rely on the following technical lemma.
\begin{lemma}\label{L4} In the regime \eqref{E11}, we have:

\noindent {\rm (i)} Let $(\tau_x: x\geq 0)$ be a inverse Gaussian subordinator with cumulant
$$\kappa(q)= \sigma^{-2}\left(\sqrt{c^2+2q\sigma^2}-c\right) = c^{-1}\int_0^{\infty} (1-\e^{-qy})\nu({\rm d}y)\,,\qquad q\geq 0\,,$$
i.e. with zero drift and L\'evy measure $c^{-1}\nu$
where  $\nu$ given in Proposition \ref{P2}. 
Then
$${\mathcal L}\left ((n^{-2}T_0,n^{-1}M_1), \P^{p(n)}_{a(n)}\right)   \ \Longrightarrow\ 
 (\tau_x, c\tau_x)\,.$$

\noindent {\rm (ii)} The behavior of the joint tail distribution of $T_0$ and $M_1$ 
under $\P^{p(n)}_{1}$ is given by
$$\lim_{n\to\infty} n \P^{p(n)}_1 \left (n^{-2}T_0> t \hbox{ or } n^{-1}M_1> m\right)   \,=\, c^{-1}
\bar \nu(\min(t,m/c))\qquad \hbox{in }L^1_{\rm loc}(\R_+\times \R_+,{\rm d}t\, {\rm d}m)\,,$$
where $\bar \nu$ denotes the tail function of the L\'evy measure $\nu$. 

\end{lemma}

\proof One could establish these limits from the explicit expressions in Proposition \ref{P1}; however a probabilistic argument based on the construction in Section 2.3 circumvents the somewhat tedious calculations.

(i) 
Recall that the fixed reproduction law $\pi\+$ has unit mean and variance $\sigma^2$. For each $n\in\N$,
consider a random walk $(S^{(n)}_k: k\in\Z_+)$ started from $S^{(n)}_0=a(n)$ 
and with step distribution that of $\xi\+-1$. By Donsker's invariance principle and Skorohod's representation, we may suppose that with probability one
$$\lim_{n\to\infty}n^{-1}S^{(n)}_{\lfloor n^2t\rfloor} =x+\sigma B_t\,,$$
where $(B_t: t\geq 0)$ is a standard Brownian motion and the convergence holds uniformly on every compact time-interval. 

For every fixed $n$, we now decompose each variable $\xi\+_i$ 
as the sum $\xi\+_i=\xi_i^{({\rm c}n)}+\xi_i^{({\rm m}n)}$ by using a Bernoulli sampling; that is conditionally on 
$\xi\+_i=\ell$, $\xi_i^{({\rm m}n)}$ has the binomial distribution with parameter $(\ell, p(n))$. 
Of course, we use independent Bernoulli sampling for the different indices $i$, so that the pairs 
$(\xi_i^{({\rm c}n)},\xi_i^{({\rm m}n)})$ are i.i.d. and have the law of $\xi$ under $\P^{p(n)}$. 
If we define
$$S^{({\rm m}n)}_k:=\xi_1^{({\rm m}n)}+ \cdots + \xi_k^{({\rm m}n)}\,,\qquad k\in\Z_+\,,$$
then $\E(\xi_1^{({\rm m}n)})= p(n)\sim c/n$ and ${\rm var}(\xi_1^{({\rm m}n)})=O(1/n)$, and it is easy to verify that with probability one
$$\lim_{n\to\infty}n^{-1}S^{({\rm m}n)}_{\lfloor n^2t\rfloor} =ct\,,$$
 uniformly on every compact time-interval. 
Hence the random walk
$$S^{({\rm c}n)}_k:=a(n)+\xi_1^{({\rm c}n)}+ \cdots + \xi_k^{({\rm c}n)}=S^{(n)}_k-S^{({\rm m}n)}_k$$
fulfills 
$$\lim_{n\to\infty}n^{-1}S^{({\rm c}n)}_{\lfloor n^2t\rfloor} =x+\sigma B_t- ct\,,$$
where again the convergence holds a.s., uniformly on every compact time-interval.

Now recall the framework of Section 2.3 and introduce
$$\varsigma^{(n)}(0):=\inf\{k\in\Z_+: S^{({\rm c}n)}_k=0\}
\quad \hbox{and}\quad \Sigma^{(n)}(0):=\sum_{i=1}^{\varsigma^{(n)}(0)}\xi^{({\rm m}n)}_i\,
=S^{({\rm m}n)}_{\varsigma^{(n)}(0)}.$$
It follows readily from the preceding observations that with probability one
$$\lim_{n\to\infty} n^{-2} \varsigma^{(n)}(0)=\tau_x \quad \hbox{and}
\lim_{n\to\infty} n^{-1} \Sigma^{(n)}(0)= c\tau_x$$
where $\tau$ denotes the process of  first passage times for a Brownian motion with drift,
$$\tau_y:=\inf\{t\geq 0: c t-\sigma B_t>y\}\,,\qquad y\geq 0.$$
It is well-known that the latter is a subordinator with cumulant $\kappa$ as given in the statement,
and the first claim is established by an appeal to Lemma \ref{L3}. 

(ii)  The branching property shows that the
law of $(T_0,M_1)$ under  $\P^{p(n)}_{a(n)}$ is that of the sum of $a(n)$ i.i.d. variables distributed as $(T_0,M_1)$ under  $\P^{p(n)}_1$. 
This observation enables to deduce (ii) from (i) by an argument similar to that we use to establish \eqref{E17}. Indeed, write
$$\bar \mu_n(t, m)=\P^{p(n)}_1\left (T_0 >t \hbox{ or } M_1>m\right)$$
for the bivariate tail distribution of the pair $(T_0,M_1)$ under $\P^{p(n)}_1$. By an elementary calculation, we have that for every $q,r>0$
\begin{eqnarray*}
& &\E^{p(n)}_{a_n}\left( \exp\left(-\frac{q}{n^2}T_0-\frac{r}{n}M_1\right)\right) \\
&=&\left(1-qr\int_0^{\infty}\int_0^{\infty} \e^{-qt}\e^{-rm}\bar \mu_n(n^2t, nm){\rm d} t\, {\rm d}m\right)^{a_n}\,.
\end{eqnarray*}
We know from (i) that this quantity converges as $n\to\infty$ towards
$$\E(\exp(-(q+cr)\tau_x))=
\exp\left(-x  \int_0^{\infty}(1-\e^{-(q+cr)y})c^{-1} \nu({\rm d}y)\right)\,,$$
so that taking logarithms, we get
\begin{eqnarray*}
& &\lim_{n\to\infty} qr a_n\int_0^{\infty}\int_0^{\infty} \e^{-qt}\e^{-rm}\bar \mu_n(n^2t, nm){\rm d} t\, {\rm d}m\\
&=&x \int_0^{\infty}(1-\e^{-(q+cr)y})c^{-1}\nu({\rm d}y)\\
&=&xqr\int_0^{\infty}\int_0^{\infty} \e^{-qt}\e^{-rm} c^{-1}\bar \nu(\min(t,m/c)){\rm d} t\, {\rm d}m\,. \end{eqnarray*}
This entails our claim.
 \QED 
 
 Proposition \ref{P2} immediately follows from Lemma \ref{L1} and Lemma \ref{L4}(i),
 so we turn our attention to the proof of Theorem \ref{T1}. 
 
{\noindent{\bf Proof of Theorem \ref{T1}:}\hskip10pt}  Recall that ${\mathcal A}_{\varnothing}=T_0$ and $d_{\varnothing}=M_1$. 
On the one hand, we know from Lemma \ref{L4}(i) that
$${\mathcal L}\left ((n^{-2}{\mathcal A}_{\varnothing},n^{-1}d_{\varnothing}), \P^{p(n)}_{a(n)}\right)   \ \Longrightarrow\ 
 (\tau_x, c\tau_x)\,.$$
On the other hand, Lemma \ref{L4}(ii) and the law of rare events for null arrays (e.g. Theorem 14.18 in \cite{Kallenberg}) entails that for any sequence of integers $b(n)$ such that $b(n)\sim bn$ for some $b>0$,
$${\mathcal L}\left ((n^{-2}T_0,n^{-1}M_1)^{(b(n)\downarrow)}, \P^{p(n)}_1\right)   \ \Longrightarrow\ 
(({\tt a}_1,c{\tt a}_1), ( {\tt a}_2, c{\tt a}_2), \ldots)\,,$$
where the notation $\gamma^{(d\downarrow)}$ has been defined just before Lemma \ref{L2} and $({\tt a}_1, {\tt a}_2, \ldots)$ stands for the sequence ranked in the decreasing order of the atoms of a Poisson measure on $(0,\infty)$ with intensity $bc^{-1}\nu$. 

Denote by $\Q_x$ the law of a tree-indexed CSBP and initial population distributed as $\tau_x$ and reproduction measure $\nu$. We now see from Lemma \ref{L2} that 
$${\mathcal L}\left (((n^{-2}{\mathcal A}_{u},n^{-1}d_{u}): |u|\leq 1), \P^{p(n)}_{a(n)}\right)   \ \Longrightarrow\ 
{\mathcal L}\left ((({\mathcal A}_{u},c{\mathcal A}_{u}): |u|\leq 1), \Q_x\right)$$
in the sense of finite dimensional convergence. Lemma \ref{L2} enables us to iterate the argument to the subsequent levels of vertices, which establishes our claim. \QED

\end{section}

\vskip 8mm

\noindent{\bf Acknowledgments}. The question of describing  the asymptotic shape of the tree of alleles for large populations with rare mutations was raised by Matthias Winkel during a lecture based on \cite{Be09} that I delivered at  the University of Oxford. I would like to thank Matthias for having stimulated this work. This work has been supported by ANR-08-BLAN-0220-01.

\end{document}